\documentstyle[12pt]{article}

\renewcommand\a{\alpha}
\newcommand\bp{\bar\p}

\newcommand\C{{\mbox{\rm\bf C\hspace{-7.3pt}}{^{_{\bf\mid}}}\hspace{4.5pt}}}
\renewcommand\d{\delta}
\newcommand\dd{\Delta}
\newcommand\D{{\cal D}}

\newcommand\g{\gamma}

\newcommand\iint{\displaystyle \int\!\!\!\int}

\newcommand\k{\kappa}
\renewcommand\l{\lambda}
\newcommand\m{\mu}

\newcommand\om{\omega}
\newcommand\op[1]{\mathop{\rm #1}\nolimits}
\newcommand\p{\partial}
\newcommand\po{$\!\!\!{\bf .}$ }
\newcommand\R{{\rm I\hspace{-2.5pt} R}}
\newcommand\st{$\!\!\!{\bf '.}$ }
\newcommand\stt{$\!\!\!{\bf ''.}$ }

\newcommand\td{\tilde}
\newcommand\te{\theta}
\newcommand\ve{\varepsilon}

\newcommand\z{\zeta}
\newcommand\Z{{\rm Z\mkern-5muZ}}

\newcommand\1{{\bf 1}}
\newcommand\qed{\phantom{\underline{y}}\hfill\hfill$\Box$}
\newcommand\bib[1]{\bibitem[#1]{#1}}
\newcommand\const{\text{const}}
\newcommand\mod{\text{mod}}
\newcommand{\text}[1]{{\mbox{\rm #1}}}
\newcommand{\dfrac}[2]{\frac{\displaystyle #1}{\displaystyle #2}}

\newtheorem{th}{Theorem}
\newtheorem{prop}{Proposition}
\newtheorem{lem}{Lemma}

\newenvironment{cor}{\trivlist \item[\hskip \labelsep{\bf Corollary.}]\it}%
{\endtrivlist}
\newenvironment{dfn}{\trivlist \item[\hskip \labelsep{\bf Definition.}]}%
{\endtrivlist}
\newenvironment{proof}{\trivlist \item[\hskip
\labelsep{{\it\underline{Proof}.\/}}]}%
{\endtrivlist}
{\endtrivlist}

\begin{document}

\title{Existence of close\\ pseudoholomorphic disks\\
for almost complex manifolds\\
and an application\\ to Kobayashi-Royden pseudonorm}
\author{B. S. Kruglikov}
\date{}
\maketitle

 \begin{abstract}
It is proved in the paper%
\footnote{Author's work was partially supported by grant INTAS 96-0713.}
 that near every pseudoholomorphic disk
on an almost complex manifold a disk of almost the same size in any
close direction passes. As an application the Kobayashi-Royden
pseudonorm for almost complex manifolds is defined and studied.
 \end{abstract}

%%%%%%%%%%%%%%%%%%%%%%%%%%%%%%%%%%%%%%%%%%%%%%%%%%%%%%%%%%%%%%%%%%%%%%%%%%%%
%0%
\section*{Introduction}

\hspace{13.5pt}
Let $(M^{2n},J)$ be an almost complex manifold, i.e.\
$J^2=-\1\in T^*M\otimes TM$. A mapping
$\Phi\:(M_1,J_1)\to(M_2,J_2)$ is called pseudoholomorphic if its
differential preserves the complex multiplication in the tangent
bundles: $\Phi_*\circ J_1=J_2\circ\Phi_*$.

Denote by $e=1\in T_0\C$ the unit vector. Let us also denote by
$D_R$ the disk in $\C$ of radius $R$, which is equipped with the
standard complex structure $J_0$. Let $v\in T_pM$,
$p=\tau_Mv$, with $\tau_M:TM\to M$ being used for the canonical
projection. We say that a disk $f:D_R\to M$ passes in the
direction $v$ if $f_*e=v$. Due to theorem III~from \cite{1}
there exists a small pseudoholomorphic disk in the direction of
an arbitrary vector $v$. We study non-small pseudoholomorphic disks
which lie in a neighborhood of a given pseudoholomorphic disk.
The main result~is

 \begin{th}\po
Let a nonconstant pseudoholomorphic disk of radius $R$
pass through a point $p$ at an almost complex manifold $(M^{2n},J)$:
 $$
f_0:(D_R,J_0)\to(M,J),\qquad(f_0)_*(0)e=v_0\ne0.
 $$
Then for every $\ve>0$ there exists a neighborhood
${\cal V}={\cal V}_\ve (v_0)$ of the vector $v_0\in TM$ such that
in the direction of each vector $v\in\cal{V}$
a pseudoholomorphic disk of radius $R-\ve$ passes:
 $$
f:(D_{R-\ve},J_0)\to(M,J),\qquad f_*(0)e=v.
 $$
 \end{th}

This theorem has an important application in the theory of invariant
metrics. In 1967 Kobayashi \cite{2} introduced a pseudodistance
on complex manifolds, which is invariant under biholomorphisms.
This gave rise to hyperbolic spaces theory [3--5].
Kobayashi pseudodistance is the maximal pseudodistance among all
pseudodistances non-increasing under holomorphic mappings, which
on the unit disk $D_1\subset\C$ coincides with the distance
$d_D$, induced by infinitesimal Poincar\'e metric in Lobachevskii
model
 $$
dl^2=\frac{dz\, d\bar z}{(1-|z|^2)^2}.
 $$
On a complex manifold $M$ the pseudodistance is defined by the formula
 $$
d_M(p,q)=\inf\sum_{k=1}^md_D(z_k,w_k),
 $$
where the infimum is taken over all holomorphic mappings
$f_k:D_1\to M$, $k=1,\dots,m$, such that $f_1(z_1)=p$,
$f_k(w_k)=f_{k+1}(z_{k+1})$ and $f_m(w_m)=q$. In the paper~\cite{6}
the Kobayashi pseudodistance was extended to the case of arbitrary
almost complex manifolds and it was shown that the basic properties
of this pseudodistance are preserved.

In 1970 Royden \cite{7} found an infinitesimal analog of the Kobayashi
pseudodistance for complex manifolds. We define the corresponding
notion in the category of almost complex manifolds and
we prove, using theorem~1, the coincidence theorem (theorem~3).
We obtain a hyperbolicity criterion (theorem~4). We also consider
the reduction procedure, which allows to define geometric invariants
of the moduli space for pseudoholomorphic curves.

%%%%%%%%%%%%%%%%%%%%%%%%%%%%%%%%%%%%%%%%%%%%%%%%%%%%%%%%%%%%%%%%%%%%%%%%%%%%
%1%
\section{\hskip-20pt .
Existence of close pseudoholomorphic disks}

%%%%%%%%%%%%%%%%%%%%%%%%%%%%%%%%%%%%%%%%%%%%%%%%%%%%%%%%%%%%%%%%%%%%%%%%%%
% 1.1 %
\subsection{\hskip-16.5pt . \hskip2pt
Reformulation of the main result}

\hspace{13.5pt}
Let us reformulate theorem~1 using the differential equation language
in appropriate coordinates. To begin with choose these coordinates
along the disk $f_0(D_R)\subset M$. Due to existence of
isothermal coordinates on surfaces~\cite{8} the disk
$f_0(D_R)$ can be defined in local complex coordinate system
$(z^1,\dots,z^n)$, which is specified in some neighborhood
of the disk, via the formulae: $|z^1|\le R$,
$z^2=\dots=z^n=0$. Moreover the disk will be pseudoholomorphic,
$J|_{\op{Im}f_0}=J_0$, and $v_0=(1,0,\dots,0)\in T_0\C\simeq\C$.

 \begin{prop}\po
In an appropriate coordinate system the vector
fields $\p_k=\p/\p z^k$ and $\bar\p_k=\p/\p\bar z^k$
at points of the disk $f_0(D_R)$ satisfy the conditions
 \begin{equation}
J\p_k=i\p_k,\ \ J\bar\p_k=-i\bar\p_k.
\label{1}
 \end{equation}
 \end{prop}

 \begin{proof}
Given equations are already satisfied on the disk for
the vector fields $\p_1$, $\bar\p_1$%
$\vphantom{\dfrac22}$.
Further at the points of
the disk we define transversal to this disk vector fields
$\p_k$, $\bar\p_k$, $k\ge2$, in such a way that all the union of
$2n$ vectors forms a basis at each point and also that
condition~(1) is satisfied. Upon constructing the needed vector
fields at the points of the disk we extend them to a neighborhood with
the help of lemma~1. Obtained structure $J$ coincides with the structure
$J_0$ on the disk and does not necessarily do so outside.
\qed
 \end{proof}

 \begin{lem}\po
Let we be given $k$ standard commuting vector fields
$v_i=\p_i$, $i=1,\dots,k$, and also $n-k$ transversal fields
$v_j$, $j=k+1,\dots,n$,
along the disk $D^k\subset\R^k\times\{0\}^{n-k}\subset\R^n$;
at each point $x\in D^k$ all the vectors $v_1,\dots,v_n$ forming a basis.
Then there exist coordinates $x^i$ in a small neighborhood
of the disk $D^k$ such that
$v_i(x)=\p_i=\p/\p x^i$, $i=1,\dots,n$, for all $x\in D^k$.
 \end{lem}

 \begin{proof}
Since the commutators of vector fields along $D^k$ are
determined by their $1$-prolongations outside $D^k$, we write the
general form for a $1$-prolongation of the vector field $v_1$:
 \begin{equation}
v_1=\sum_{r=1}^n\Bigl(\d_1^r+\sum_{s=k+1}^n
x^s \phi_s^r(x_1,\dots,x_k)\Bigr)\p_r\ \mod\, \m^2\D,
\label{2}
 \end{equation}
where $\mu^2\cal{D}$ is the submodule of the module of vector fields,
consisting of the vector fields vanishing on the submanifold
$D^k\subset\R^n$ to the second order. If on the disk $D^k$
the decomposition of the additional vector fields is written as
 \begin{equation}
v_j=\sum_{s=1}^n a_j^s(x_1,\dots,x_k)\p_s,\ j=k+1,\dots,n,
\label{3}
 \end{equation}
then the equations $[v_1,v_j]=0$ with $x^{k+1}=\dots=x^n$ have the
following form
 $$
\sum_{r=1}^n(\p_1a_j^r(x^1,\dots,x^k))\p_r=
\sum_{s=k+1}^na_j^s(x^1,\dots,x^k)
\sum_{r=1}^n\phi_s^r(x^1,\dots,x^k)\p_r.
 $$

This system decomposes (by $r$) on $n$ determinate systems of
$n-k$ linear equations with $n-k$ unknowns. The matrix
$(a_j^s)_{k+1\le j,s\le n}$ of each system is nondegenerate, hence the
system possesses a solution.

Thus the field $v_1$ is constructed. Let us rectify it:
$v_1=\p/\p x^1$. We can assume that on the disk
$D^k\subset\{x^{k+1}=\dots=x^n=0\}$ the tangent vector fields have
the original form $v_i=\p_i$. In new coordinates the
coefficients of the decomposition~(3) do not depend on $x^1$.
So one can search for the prolongation of the field $v_2$ in
the form similar to~(2), but with no dependence on $x^1$.
Continuing the process we get some coordinates $x^1,\dots,x^n$,
in which the disk $D^k$ belongs to the subspace
$\{x^{k+1}=\dots=x^n=0\}$ and such that on this disk
 \begin{equation}
v_i=\dfrac\p{\p x^i},\,i=1,\dots,k,\
v_j=\sum_{s=1}^n a_j^s\p_s,\ a_j^s=\const,\,j=k+1,\dots,n.
\label{4}
 \end{equation}
Now we prolong the vector fields to a neighborhood via the
formula~(4).
 \qed
 \end{proof}

Thus we introduce complex coordinates $z^k=x^k+iy^k$ in a neighborhood
of the disk $f_0(D^k)$. Now our manifold, being contracted,
has the form
 \begin{equation}
M_0=\{|z_1|\le R,\ |z_k|\le R_1,\,k\ge2\}\simeq D_R\times(D_{R_1})^{n-1}
\subset\C^n,\ \
R_1\ll R,
 \label{5}
 \end{equation}
and the structure $J$ at points of the disk
$D_R=\{(x_1,y_1,0,\dots,0,0)\}$ has the form
 \begin{equation}
J\dfrac\p{\p x^k}=\dfrac\p{\p y^k},
J\dfrac\p{\p y^k}=-\dfrac\p{\p x^k}.
 \label{6}
 \end{equation}

Writing down Cauchy-Riemann equations $f_*\circ J_0=J\circ f_*$ on
the mapping of the disk $f:(D_{R-\ve},J_0)\to(M,J)$
(similar to sec.~3.3 from \cite1) and using the rectifying conditions~(6),
we get an equivalent formulation of the main statement
($\p,\bar\p$ are considered in the coordinates $z^k$):

 \addtocounter{th}{-1}
 \begin{th}\st
Let $n^2$ functions
$a_{\bar m}^i$ {\rm(}$i,m=1,\dots,n${\rm)} of the
class $C^{k+\l}$, $k\in\Z_+$, $\l\in(0,1)$ be given
on a manifold $M_0$ of the form\/~{\rm(5)}.
Let $\ve\in(0,R)$ be an arbitrary small real number. If
$a^i_{\bar m}(z)=0$ for all points $z\in D_R\times\{0\}^{n-1}\subset M_0$,
then the equation
 $$
\bar\p z^i+\sum_{m=1}^n a_{\bar m}^i(z)\bar\p\bar z^m=0,\quad
z^i(0)=p^i,\quad\p z^i(0)=u^i,\quad i=1,\dots,n,
 $$
has a solution $z^i=z^i(\zeta)\in C^{k+1+\l}(D_{R-\ve};M_0)$
subject to the restriction that the neighborhood
${\cal V}={\cal V}(v_0)\ni v$ of the vector
$v_0=(1,0,\dots,0)\in T_0M_0$ is chosen sufficiently small.
Here $v=(p,u)$, $p=\tau_Mv\in M_0$, $u\in T_pM_0$.
 \end{th}

%%%%%%%%%%%%%%%%%%%%%%%%%%%%%%%%%%%%%%%%%%%%%%%%%%%%%%%%%%%%%%%%%%%%%%%%%%
% 1.2 %
\subsection{\hskip-16.5pt . \hskip2pt
Covering of the neighborhood by disks and another reformulation}
 \addtocounter{th}{-1}
 \begin{th}\stt
One can set $p=0$ in the formulation of
theorem $1'$. Thus in the chosen coordinates
the equation on the pseudoholomorphic disk sought for takes
the following form\/{\rm:}
 \begin{equation}
\left\{
\begin{array}{rcl}
\bar\p z^1&=&-\sum_{m=1}^n a_{\bar m}^1(z)\bar\p\bar z^m,\\
\bar\p z^I&=&-\sum_{m=1}^n a_{\bar m}^I(z)\bar\p\bar z^m.
\end{array}
\right.
 \label{7}
 \end{equation}
 $$
z^1(0)=0,z^I(0)=0,\ \big(\p z^1(0),\p z^I(0)\big)=(u^1,\dots,u^n),
 $$
where the multiindex $I$ stands for $(2,\dots,n)$.
 \end{th}

Next statement shows equivalence of theorems $1'$ and $1''$.

 \begin{prop}\po
For every pseudoholomorphic disk
$f_0:D_R\to M_0$ and every $\ve>0$ a small neighborhood of the image
$f_0(D_{R-\ve})$ can be covered by the images of close
pseudoholomorphic disks $f$ of radii $R-\d$, where $\d<\ve$\rm:
${\cal O}(\op{Im}f_0(D_{R-\ve}))\subseteq\bigcup_f\op{Im}f(D_{R-\d})$.
 \end{prop}

 \begin{proof}
Perturb the almost complex structure $J$ in a
neighborhood of the disk $f_0(D_R)$ so that it coincides with
the standard integrable structure $J_0$ near the boundary of this
neighborhood: for every $\ve>0$ there exists such an almost
complex structure $\td{J}$ that $\td{J}=J$ in a small neighborhood
of the disk $f_0(D_{R-\ve/2})=D_{R-\ve/2}\times\{0\}^{n-1}\subset M_0$
and $\td{J}=J_0$ in a neighborhood of the boundary of the manifold
$\td{M}=D_R\times(D_\d)^{n-1}\subset M_0\subset\C^n$.
Further one can suppose that
$\td{M}\subset\hat{M}=S^2_R\times(S^2_{2R})^{n-1}\simeq(S^2)^n$,
where $\hat{M}$ can be equipped with an almost complex structure
$\hat{J}$, which coincides with $\td{J}$ in $\td{M}$ and which equals
the standard integrable structure $J_0$ in the complement.
Let us supply the manifold $\hat{M}$ with the symplectic structure
$\om=\om_0^{(1)}\oplus\om_0^{(2)}\oplus\dots\oplus\om_0^{(2)}$, where
$\om_0$ is the standard volume form, and also
$\om_0^{(1)}(S^2_R)=\pi R^2$ and $\om_0^{(2)}(S^2_{2R})=4\pi R^2$.
Decreasing if necessary the size of the neighborhood of the
disk $f_0(D_{R-\ve/2})$ we can suppose that the almost complex
structure $\hat{J}$ is tamed by the symplectic structure
$\om$, i.e.\ $\om(\xi,\hat{J}\xi)>0$ for $\xi\ne0$.

Denote by $A\in H_2(\hat{M};\Z)$ the homology class of
the sphere $S^2_R\times\{*\}^{n-1}\subset\hat{M}$.
The disk $f_0(D_{R-\ve/2})$ can be extended
to the entire rational pseudoholomorphic curve
$u_0:S^2\to\hat{M}$, which lies in the class $A$.
Let us consider the space ${\cal M}(A,\hat{J}\,)$
of entire pseudoholomorphic curves $u:S^2\to\hat{M}$ of the class~$A$.
Since the class $A$ cannot be decomposed into a sum of homology
classes $\sum_{i=1}^nA_i$, $n\ge2$, with $\om(A_i)>0$, then
Gromov compactness theorem (\cite{9}~sec.~1.5.B or \cite{10}~Sec.~4.3.2)
implies the compactness of the space ${\cal M}(A,\hat{J}\,)/G$, where
$G\simeq PSL_2$ is the complex automorphisms group of the sphere
$(S^2,J_0)$, $\op{dim}G=6$. Moreover for almost complex structure
$\hat{J}$ of general position the space ${\cal M}(A,\hat{J}\,)$
is a smooth manifold of dimension $2n+4$
\cite{9}~sec.~2.1--2.2; \cite{10}~sec.~3.1.2.
Consider the space of nonparametrized pseudoholomorphic curves
${\cal W}\,(A,\hat{J}\,)={\cal M}(A,\hat{J}\,)\times_G S^2$.
This space is a compact manifold of dimension $2n$. Let us consider
the evaluation map $e:{\cal W}\,(A,\hat{J}\,)\to\hat{M}$, which is
defined by the formula $e(u,z)=u(z)$ for $z\in S^2$,
$u\in{\cal M}(A,\hat{J}\,)$. We suppose the group $G$ acts on
${\cal M}\times S^2$ by conjugation,
$\phi(u,z)=(u\circ\phi^{-1},\phi(z))$, $\phi\in G$, whence the
correctness of the definition for $e$. Since $A$-curves foliate
the manifold $\hat{M}$ outside a small neighborhood of the image
$u_0(S^2)$ (because there $\hat{J}=J_0$), the map
$e$ has degree $1$, $\op{deg}e=1$. Therefore through every point,
close to the curve $u_0(S^2)$, some pseudoholomorphic curve
$u(S^2)$ passes, which is homologous to the curve $u_0(S^2)$.

To eliminate the general position condition for $\hat{J}$
we take a sequence $\hat{J}_k$ of the general position almost
complex structures, which tends to $\hat{J}$ in $C^\infty$-topology,
and use the compactness theorem \cite{10}B.4.2.
Intersecting the obtained set of pseudoholomorphic spheres
$u(S^2)$ with a small neighborhood ${\cal O}$ of the disk
$f(D_{R-\ve/2})$, we get the desired set of the disks
$f(D)$ in a neighborhood $({\cal O},J)$.
Smoothness of these disks follows from the standard elliptic
regularity \cite{1} sec.~5.4, 4.3; \cite{10} sec.~B.4.1.
 \qed
 \end{proof}

%%%%%%%%%%%%%%%%%%%%%%%%%%%%%%%%%%%%%%%%%%%%%%%%%%%%%%%%%%%%%%%%%%%%%%%%%%
% 1.3 %
\subsection{\hskip-16.5pt . \hskip2pt
Spaces, norms and estimates}

\hspace{13.5pt}
Define the $\l$-H\"older norm of complex-valued functions on the
disk $D_R$ of radius $R$ by the formula
$\|f\|=|f|+(2R)^\l H_\l[f]$, $\l\in(0,1)$,
where
 $$
H_\l[f]=\sup_{w\ne0}\bigg|\frac{f(z+w)-f(z)}{w^\l}\bigg|,\qquad
|f|=\sup|f(z)|.
 $$
The space $C^\l(D_R,M_0)$ of $\l$-H\"older maps consists of
all maps $f:D_R\to M_0$, the components of which have finite
$\l$-norms, $\|f^i\|<\infty$. The space
$C^{k+\l}(D_R,M_0)$, $k\in\Z_+$,
of $(k+\l)$-H\"older maps consists of all maps, the partial
derivatives of which up to the $k$-th order inclusive belong to
$C^\l$.

Let us also introduce the space $B=C_\bullet^{k+\l}(D_R,M_0)$
consisting of all maps $f\in C^{k+\l}$, $f(0)=0$, with the norm
$\|f\|=\sum_1^n\|f^i\|$. Note that the space $B$ can be also
supplied with the norm
 $$
\|f\|'=\max\>\{\|\p f\|,\|\bar\p f\|\}.
 $$

 \begin{prop}\po
The spaces $(C^\l,\|\cdot\|)$ and
$(C_\bullet^{1+\l},\|\cdot\|')$ are Banach.
 \end{prop}

The second statement follows from the first and the estimate
\cite{1}~7.1.c--7.1.e
 \begin{equation}
\|f\|\le6R\|f\|'.
 \label{8}
 \end{equation}

Consider the Cauchy operators
 $$
Sf(w)=\frac1{2\pi i}\oint_{\p D_R}\frac{f(\zeta)}{\zeta-w}\,d\zeta,\qquad
Tf(w)=\frac1{2\pi i}\iint_{D_R}\frac{f(\zeta)}{\zeta-w}\,d\zeta\wedge
d\bar\zeta.
 $$

Recall the basic properties of these operators \cite{1}~6.1--6.2:
 \begin{equation}
f\in C^\l(D) \Rightarrow Tf\in C^{1+\l}(D),\ \bar\p Tf=f,
 \label{9}
 \end{equation}
 \begin{equation}
f\in C^\l(D) \Rightarrow Sf\in C^\l(\op{Int}D),\ \bar\p Sf=0,\ STf=0,
 \label{10}
 \end{equation}
 \begin{equation}
f=Sf+T\bar\p f\quad
\text{ (Cauchy-Green-Pomp\'eiu formula) },
 \label{11}
 \end{equation}
 \begin{equation}
\|Tf\|'\le c_1\|f\|,\quad
\|Sf\|\le c_2\|f\|,
 \label{12}
 \end{equation}

Consider also the operators
$T_kf(w)=Tf(w)-\sum_{s=0}^k\frac1{s!}\p^s Tf(0)w^s$.

 \begin{lem}\po
For the points $w\in\op{Int}D$ the following formula holds
 $$
T_kf(w)=\frac{w^{k+1}}{2\pi i}\iint_{D_R}
\frac{f(\zeta)}{(\zeta-w)\zeta^{k+1}}\,d\zeta\wedge d\bar\zeta.
 $$
 \end{lem}

 \begin{lem}\po
The operator $T_\infty=\lim_{k\to\infty}T_k$ is defined
for functions $f\in C^\l(D_R)$, and moreover
$T_\infty f\in C^{1+\l}(D_{R-\ve})$.
 \end{lem}

 \begin{lem}\po
 $
T_k(w^l\bar w^m)=
 \left[\begin{array}{l}
\dfrac{w^l\bar w^{m+1}}{m+1}, \quad l<k+m+2, \\
\dfrac{w^l\bar w^{m+1}}{m+1}-\dfrac{R^{2(m+1)}}{m+1}w^{l-m-1},
\quad l\ge k+m+2.
 \end{array}\right.
 $
\label{L4}
 \end{lem}

 \begin{cor}
$T_\infty(w^l\bar w^m)=w^l\bar w^{m+1}\!/(m+1)$.
 \end{cor}

Thus the operator $T_\infty$ represent the integration by
$\bar\zeta$ of the polynomials on $D_R$. Let us also
besides the space $B$ consider its closed subset
$B_\d=\{f=(f_1,\dots,f_n)\in
B,\,|f_1-\zeta|\le\d,\,|f_k|\le\d,\,k\ge2\}$.
We will seek a solution $f$ of the Cauchy-Riemann equation~(7)
in the space $B_\d$ for a small neighborhood $\cal{V}$ of the
vector $v_0$.

%%%%%%%%%%%%%%%%%%%%%%%%%%%%%%%%%%%%%%%%%%%%%%%%%%%%%%%%%%%%%%%%%%%%%%%%%%
% 1.4 %
\subsection{\hskip-16.5pt . \hskip2pt
Proof of theorem 1$''$}

\hspace{13.5pt}
{\it Idea of the proof.\/} Equation~(7) was solved in the
paper~\cite{1}, theorem III, where the velocity vector $v$ was fixed
and the radius $R\ll1$ of the disk was supposed small.
For this the linearization of almost complex structure at the point
was considered. Because of the proximity of equations on
pseudoholomorphic curves for the given almost complex structure $J$
and for the linearized one $J_0$ the following map was contractible:
 \begin{equation}
\Phi: B\to B,\ \
(\Phi f)^i(\z) = v^i\z +
T_1\left(-\sum_m a^i_{\bar m}(f) \bar\p {\bar f}^m\right)(\z).
 \label{13}
 \end{equation}
In our situation radius of the disk is not small, therefore the
word for word carrying over the arguments from~\cite{1} is possible
only if the structure $J$ differs from $J_0$ on the disk $f_0$ by
a second order smallness quantity, i.e.\ if the functions
$a_{\bar1}^i$ on the disk $f_0$ as well as their derivatives
vanish. In general situation it is not the case, so we linearize
the almost complex structure $J$ along the disk~$f_0$. Here
the linearization is parametrized by the coordinate $z^1=\zeta$
along this disk. Solutions of complex linear equation behave
similarly to solutions of the real equation $\dot x=Ax+B$:
for non-small values of the parameter it is false that
$e^{At}\approx1$, so the terms of the series
$e^{At}=\sum_{s=0}^\infty(At)^s\!/s!$ are not absolute decreasing,
but this property becomes true beginning with some number
$s\ge s_0$. Thus finite sums of the series for the exponent does
not form a contracting sequence, yet to achieve this one should
consider the sums beginning with some big number.

Let us turn to the proof.
Similarly to~\cite1, starting with formulae~(11), (10),
we seek a solution of equation~(7) in the form~(13), but we
replace the space $B$ by~$B_\d$. In fact, as noted above,
we should change the definition of the operator
$\Phi$ to improve the convergence. Let us consider the automorphism
of the space~$\C^n$, which comes from the contraction of the
space $M_0$,
 \begin{equation}
z^1\to z^1,\ z^I\to\dfrac {z^I}N,\ N\gg1.
 \label{14}
 \end{equation}
Since $a^i_{\bar m}=0$ along the disk
$D_R\times\{0\}^{n-1}\subset M_0$,
the function $a^1_{\bar1}$ becomes small and the functions
$a^I_{\bar1}$ become very close to their linearizations by the
variables $z^I$ in the norm $\|\cdot\|'$
for large $N$ in equation (7).

Consequently the first equation of~(7), considered as one
$z^I$-parametric equation, can be solved by the iteration method,
when we use formula~(13) for complex dimension $1$ and change
$B$ to $B_\d$. For small $\d$ and big $N$ in~(14) the estimates
from \cite{1} sec.~5.2 yield the contractibility
of the iteration procedure in the norm $\|\cdot\|'$. This
iteration procedure will be denoted by $z^1\mapsto\Psi^1(z^1,z^I)$.

To consider the second equation of~(7) let us linearize the
functions used in it by $z^I$:
 \begin{equation}
a^I_{\bar 1}(z)=\sum\limits_{m\ge 2}
\left(a^I_{\bar 1;m}(z^1) z^m +
a^I_{\bar 1;\bar m}(z^1) \bar z^m\right)
+ \hat a^I_{\bar 1}(z).
 \label{15}
 \end{equation}
In this formula the functions $\hat{a}^I_{\bar1}(z)$ have the second
order of smallness along the disk $D_R\subset M_0$. Let us also set
$\hat{a}^I_{\bar m}(z)=a^I_{\bar m}(z)$ when $m\ne1$.

According to Weierstrass theorem the coefficients at linear by
$z^I$ terms in~(15) are approximated by polynomials depending on
$z^1$, $\bar z^1$ in the norm $|\cdot|$ on $D_R$:
 $$
a^I_{\bar 1;m}(z^1)= p^I_m(z^1,\bar z^1) + \a^I_m(z^1),\
a^I_{\bar 1;\bar m}(z^1)= p^I_{\bar m}(z^1,\bar z^1) + \a^I_{\bar
m}(z^1),\
|\a^I_m|,|\a^I_{\bar m}|<\ve.
 $$

Let $A^I(\zeta,z^I)=\sum_{m\ge2}(p_m^I(\zeta,\bar\zeta)z^m+ p_{\bar
m}^I(\zeta,\bar\zeta)\bar z^m)$,
$A^I_\d(\zeta,z)=A^I(\zeta,z^I)-A^I(z^1,z^I)$.
Then the second equation of~(7) can be written in the form
 \begin{equation}
\bp z^I(\z)=-A^I(\z,z^I)+ U^I(z(\z)),
 \label{16}
 \end{equation}
where the summands of the remainder
$U^I=A^I_\d+U_1^I+U_2^I+U_3^I$ have the form
 $$
U_1^I(z)=A^I(z^1,z^I) (1-\bp \bar z^1),\ \
U_2^I(z)=-\sum_m \hat a^I_{\bar m}\bp\bar z^m,
 $$
 $$
U_3^I(z)=-\sum_{m\ge 2}(\a^I_m(z^1)z^m+\a^I_{\bar m}(z^1)\bar z^m)
\bp\bar z^1.
 $$

We approximate equation~(16) by the following equation with linear by
$z^I$ right hand size and polynomial by $\zeta$ coefficients:
 \begin{equation}
\bp z^I(\z)=-A^I(\z,z^I).
 \label{17}
 \end{equation}

A solution of this equation can be constructed as the limit of the
iteration procedure
 \begin{equation}
z^I_{(k+1)}= v^I\z - T_\infty[A^I(\z,z^I_{(k)})].
 \label{18}
 \end{equation}

By the corollary of lemma~4 the iteration of integration by means
of the operator $T_\infty$ has the form
$T_\infty^k(w^l\bar w^m)=w^l m!\,\bar w^{m+k}\!/(m+k)!$,
which implies that the iteration process~(18) converges under
any initial condition $z^I_{(0)}$ to a solution of equation~(17),
and moreover the convergence is exponential.
In particular, beginning with some number $k$,
the sequence $z^I_{(k)}$ is contractible. And what is more
there exist constants $C$ and $\mu$, depending only on
almost complex structure $J$ (i.e.\ on coefficients $a^i_{\bar m}$),
such that for every $k\ge1$ and polynomial $p(\zeta,\bar\zeta)$
the following inequality holds:
 \begin{equation}
\|T_\infty^k[A^I(\z,p)]\|'\le C e^{\mu R} \|p\|.
 \label{19}
 \end{equation}

We now define the iteration procedure to compute $z^I(\zeta)$.
Let the iterative term $z^I_{[r]}$ be already constructed.
Additionally in virtue of the previous step the given term
is equal to the sum of a polynomial
$P^I_{[r]}(\zeta,\bar\zeta)$ and a function
$\te^I_{[r]}(\zeta)\in C^{1+\l}$.
Represent the last function by Weierstrass theorem
as the sum of a polynomial (by $\zeta$, $\bar\zeta$)
and an error:
$\te^I_{[r]}(\zeta)=Q^I_{[r]}(\zeta)+q^I_{[r]}(\zeta)$,
$|q^I_{[r]}|\le\nu|\te^I_{[r]}|$.
Define the next term by the formula
 $$
z^I_{[r+1]}(\zeta)=v^I\zeta-T_\infty[A^I(P^I_{[r]})]
-T_\infty^{k_r}[A^I(Q^I_{[r]})]-T_1[A^I(q^I_{[r]})]+T_1[U^I(z_{[r]})].
 $$
Here $A^I=A^I(\zeta,\cdot)$ and $k_r$ is such a number that
beginning with number $k_r$ the sequence $T^k_\infty[A^I(Q^I_{[r]})]$
contracts with the coefficient $\ve_r$. In addition (cf.~(8))
the following estimates for the additional terms take place:
 $$
\|A^I_\d(\z,z')-A^I_\d(\z,z'')\|\le c_3\d\|z'-z''\|',\
\|U^I_1(z')-U^I_1(z'')\|\le c_3\d\|z'-z''\|',
 $$
 $$
\|U^I_2(z')-U^I_2(z'')\|\le c_3\d\|z'-z''\|',\
\|U^I_3(z')-U^I_3(z'')\|\le c_4\ve\|z'-z''\|'.
 $$
Taking inequality~(12) into account we conclude that for small
$\d$, $\ve$, $\ve_r$ and $\nu$ the sequence $z^I_{[r]}$ is
contractible:
$\|z^I_{[r+1]}-z^I_{[r]}\|'\le(1-\k)\|z^I_{[r]}-z^I_{[r-1]}\|'$ for
some $\k<1$ independent of $r$. Therefore, taking into consideration
the iteration by $\Psi^1$ for the variable $z^1$, we get
a convergent in $C^{1+\l}$ sequence, the limit of which has to be
the desired solution. Actually, set
$z^1_{[r+1]}=\Psi^1(z^1_{[r]},z^I_{[r]})$, taking as parameter
$z^I$ the iterative term $z^I_{[r]}$. In what follows in
determination of the term $z^I_{[r+1]}$ we assume
$z^1=z^1_{[r]}$. Thus we obtain the sequence $z_{[r]}$.

Due to the estimates considered and inequality~(19) the terms and
the limit of the sequence $z_{[r]}$ differ from its initial term
$z_{[0]}=v\zeta$ less than exponentially by $R$ with respect to
$|v-v_0|$ in the norm $\|\cdot\|'$. Therefore for small
$|v-v_0|\ll1$ all the terms and the limit of the iterative
sequence lie in $B_\d$. Hence the sequence converges in~$B_\d$.
Now it is easily seen that the limit of the sequence
$z_{[r]}$ is a solution of the equation~(7). When the coefficients
have smoothness $a^i_{\bar m}\in C^{k+\l}(M_0)$, then the
obtained solution, which is of smoothness $C^{1+\l}$, will be
actually of higher smoothness class $C^{k+1+\l}$. This follows
from the standard elliptic regularity methods for our equation
\cite{1} 5.4, 4.3, \cite{10}~B.4.1. When
$a^i_{\bar m}\in C^\infty(M_0)$
we get a smooth solution of the Cauchy-Riemann equation
$z(\zeta)\in C^\infty(D_{R-\ve};M_0)$. \qed

%%%%%%%%%%%%%%%%%%%%%%%%%%%%%%%%%%%%%%%%%%%%%%%%%%%%%%%%%%%%%%%%%%%%%%%%%%
% 1.5 %
\subsection{\hskip-16.5pt . \hskip2pt
Jet spaces and connection with $h$-principle}

\hspace{13.5pt}
Let us call the foliation by pseudoholomorphic disks any
embedding (immersion) $\Phi:D_R\times N^{2n-2}\to M$
such that all the mappings $\Phi|_{D_R\times\{x\}}$ are
pseudoholomorphic and the image of the map $\Phi$ covers
the entire manifold $M$. The construction of proposition~2
together with the positivity of intersections in dimension $4$
(\cite{9} 2.1.C$_2$; \cite{11}~1.1) imply

 \begin{prop}\po
Let $(M,J)$ be a four-dimensional almost
complex manifold. For every embedded\/ {\rm(}immersed\/{\rm)}
pseudoholomorphic disk $f:D_R\to M$ and every $\ve>0$
small neighborhood of the image $f(D_{R-\ve})$ allows the
foliation by pseudoholomorphic disks.
\label{P}
 \end{prop}

Let us consider the manifold of pseudoholomorphic jets
${\cal J}^{1}_{PH}(D_R;M)$ of the mappings $u:D_R\to M$.
Its points are triples $(\zeta,z,\Phi)$, where $\zeta\in D_R$,
$z\in M$, and $\Phi:(T_\zeta D_R,J_0)\to (T_zM,J(z))$ is a
complex linear mapping. It was shown in the paper~\cite{13} that
the manifold ${\cal J}^{1}_{PH}$ possesses a canonical
almost complex structure $J_{[1]}$, which is equal to
$J_0\oplus J\oplus J$ regarding the induced by some minimal
connection decomposition
$T_p{\cal J}^{1}_{PH}=T_\zeta D_R\oplus T_zM\oplus T_p{\cal F}$,
where $\cal{F}$ is the fiber of the natural projection
$\tau:{\cal J}^{1}_{PH}(D_R,M)\to D_R\times M$.
The canonical projection $\pi:{\cal J}^{1}_{PH}(D_R;M)\to M$
is pseudoholomorphic and any pseudoholomorphic mapping
$f:D_R\to M$ lifts canonically to the pseudoholomorphic mapping
$j^1f:D_R\to{\cal J}^{1}_{PH}(D_R;M)$,
$j^1f(\zeta)=(\zeta,f(\zeta),d_\zeta f)$.

We define the structure $J_{[1]}$ in a different way (cf.~\cite{14},
remark~1). If $p=(\zeta,z,\Phi)\in{\cal J}^{1}_{PH}$, we can
assume that the mapping $\Phi$ is the differential at the point
$\zeta$ of some small pseudoholomorphic disk $u:D_\ve\to M$.
Denote by $p^{(2)}$ the 2-jet of the disk $u$ at the point
$\zeta\in D_\ve\subset D_R$. Consider the map
$j^1u:D_\ve\to{\cal J}^{1}_{PH}$. The tangent space at the point
$p$ depends only on the value $p^{(2)}$. Denote this tangent space by
$L_{p^{(2)}}$.

Consider the natural projection $\rho:{\cal J}^{1}_{PH}\to D_R$ with
the fiber $\cal{H}$. We have
$T_p{\cal J}^{1}_{PH}=L_{p^{(2)}}\oplus T_p\cal{H}$,
both summand being naturally equipped with complex structures.
Set $J_{[1]}=J_0\oplus J$. This structure does not depend on the choice
of $p^{(2)}$, i.e.\ it is defined canonically.

 \begin{dfn}
Let us call a pseudoholomorphic disk
$g:D_R\to{\cal J}^{1}_{PH}(D_R,M)$ {\it holonomic\/},
if the mapping $g$ is the 1-jet lifting of some pseudoholomorphic disk
from $D_R$ to $M$: $g=j^1f$.
 \end{dfn}

Proposition~2 applied to a holonomic disk
$g=j^1f:D_R\to{\cal J}^{1}_{PH}$ yields existence of a
pseudoholomorphic disk $g'$ through each point arbitrary close to the
image of the disk $g$, which however needs not be a holonomic disk,
$g'\ne j^1(\pi\circ g')$. In this sense theorem~1 provides a more
strong statement. Actually, closeness of initial points of the disks
$g=j^1f$ and $g'=j^1f'$ in ${\cal J}^{1}_{PH}(D_{R-\ve};M)$ means
closeness of initial points and initial directions of the maps
$f$ and $f'$ in $TM$. Thus theorem~1 implies existence of
$C^1$-close disk $f'$, and we can set $g'=j^1f'$. Thus we proved

 \begin{th}\po
Through every point, which is close to the image of embedded\/
{\rm(}immersed\/{\rm)} holonomic pseudoholomorphic disk
$g:D_R\to{\cal J}^{1}_{PH}(D_R;M)$, an embedded\/
{\rm(}immersed\/{\rm)} holonomic pseudoholomorphic disk
$g':D_{R-\ve}\to{\cal J}^{1}_{PH}(D_{R-\ve};M)$ passes. \qed
 \end{th}

In other words, proposition~2 remains also valid in the holonomic
situation. The statement just proved is a particular case of
the so-called $h$-principle \cite{15}. It is also interesting to
get the holonomic version of proposition~4.

%%%%%%%%%%%%%%%%%%%%%%%%%%%%%%%%%%%%%%%%%%%%%%%%%%%%%%%%%%%%%%%%%%%%%%%%%%%%
%2%
\section{\hskip-20pt . \hskip2pt
Kobayashi-Royden pseudonorm}

%%%%%%%%%%%%%%%%%%%%%%%%%%%%%%%%%%%%%%%%%%%%%%%%%%%%%%%%%%%%%%%%%%%%%%%%%%
% 2.1 %
\subsection{\hskip-16.5pt . \hskip2pt
Definition of the pseudonorm and its main properties}

\hspace{13.5pt}
Let us consider the set ${\cal R}(v)=\bigcup_{r>0}{\cal R}_r(v)$,
where ${\cal R}_r(v)$ for $r\in\R_+$ consists of pseudoholomorphic
mappings $f:D_1\to M$, such that $f_*(0)e=rv$.

 \begin{dfn}
Let us call the Kobayashi-Royden pseudonorm on an almost
complex manifold $M$ the function on the tangent bundle $TM$,
which is defined by the formula
 $$
F_M(v)=\inf_{{\cal R}(v)}\frac1r.
 $$
 \end{dfn}

According to theorem III from~\cite{1} the set ${\cal R}_r(v)$ is
nonempty for small $r$, so the definition is correct. We call the
function $F_M$ pseudonorm since it is nonnegative and homogeneous of
degree one: $F_M(tv)=|t|F_M(v)$. However $F_M$ can vanish in some
directions and the triangle inequality does not hold. The
next statement follows from the very definition.

 \begin{prop}\po
Given any vector $v\in TM_1$ and any
pseudoholomorphic mapping $f:(M_1,J_1)\to(M_2,J_2)$ we have
$$
F_{M_2}(f_*v)\le F_{M_1}(v).
$$
 \end{prop}

Let us fix some norm $|\cdot|$ on $TM$.

 \begin{prop}\po
\label{A}
{\rm(i)} There exists a constant $C_K$ for every
compact $K\subset M$ such that each vector $v\in TM$ with
$\tau_M v\in K$ satisfies
 $$
F_M(v)\le C_K|v|.
 $$
{\rm(ii)} Let $M$ be a compact manifold\/ {\rm(}with possible
boundary\/{\rm)} equipped with an almost complex structure $J$,
which is tamed by an exact symplectic form $\om=d\a$,
$\om(\xi,J\xi)>0$ for $\xi\ne0$. Then there exists such a constant
$c_M>0$, that for all $v\in TM$
$$
F_M(v)\ge c_M|v|.
$$
 \end{prop}

 \begin{proof}
For a small neighborhood $U$ of the point $p\in M$
the estimates of sec.~5.2a of the paper \cite{1} imply existence
of a number $\ve>0$, dependent only on the almost complex structure
$J$ and the neighborhood $U$, such that for every
$q\in U$, $v\in T_qM$, $|v|=1$, and $r\in(0,\ve)$
there exists a pseudoholomorphic disk $f:D_1\to M$ such that
$f(0)=q$, $f_*(0)e=rv$. Setting $C_U=1/\ve$ we have
$F_M(v)\le C_U|v|$ for all (now not necessarily unit) vectors $v$
for which $\tau_Mv\in U$. Since a compact set can be covered by
a finite number of neighborhoods $U$, the first statement of the
proposition is proved. The second part is a reformulation of the
nonlinear Schwarz lemma \cite{9}~1.3.A: if an almost complex
structure $J$ on a compact manifold is tamed by an exact symplectic
structure $\om$, then the derivative at zero of any pseudoholomorphic
disk $f:D_1\to M$, passing through a fixed point at the manifold,
is bounded by a non-depending on the disk constant: $|f_*(0)e|<C$.
 \qed
 \end{proof}

 \begin{prop}\po
The function $F_M$ is upper semicontinuous.
\label{B}
 \end{prop}

 \begin{proof}
The inequality $\overline{\lim\limits_{v\to v_0}}F_M(v)\le F_M(v_0)$
is equivalent to the statement of theorem~1 because $F_M(v)=\inf(1/R)$,
where the lower bound is considered over all mappings $f:D_R\to M$,
such that $f_*(0)e=v$.\qed
 \end{proof}

%%%%%%%%%%%%%%%%%%%%%%%%%%%%%%%%%%%%%%%%%%%%%%%%%%%%%%%%%%%%%%%%%%%%%%%%%%%%
%2.2%
\subsection{\hskip-16.5pt . \hskip2pt Coincidence theorem}

\hspace{13.5pt}
Define a function ${\bar d}_M:M\times M\to\R$ by the formula
 $$
{\bar d}_M(p,q)=\inf_\g\int_0^1F_M(\dot\g(t))\,dt,
 $$
where the lower bound is taken over all piecewise smooth paths $\g$
from the point $p$ to~$q$. Propositions~6(i) and~7 imply
correctness of the definition and

 \begin{prop}\po
\label{Pr8}
The function ${\bar d}_M$ is pseudodistance. \qed
 \end{prop}

 \begin{th}\po
Introduced pseudodistance coincides with the Kobayashi
pseudodistance, $d_M={\bar d}_M$.
 \end{th}

 \begin{proof}
The inequality ${\bar d}_M\le d_M$ is evident because
$F_M(v)=\inf|\xi|$, where the lower bound is taken over all
pseudoholomorphic mappings $f:D_1\to M$, $f_*\xi=v$,
and the norm is count with respect to the Poincar\'e metric.
Let us prove the reverse. We follow the Royden's proof \cite{7}.

Let $\g$ be a smooth curve from a point $p$ to a point $q$ such that
$\int_\g F_M<{\bar d}_M(p,q)+\ve$. Due to upper semicontinuity
there exists a continuous on $[0,1]$ function $h$, such that
$h(t)>F_M(\dot\g(t))$ and
 $$
\int_0^1h(t)\,dt<{\bar d}_M(p,q)+\ve,
 $$
i.e.\ for sufficiently dense partition $0=t_0<t_1<\dots<t_k=1$ we have
 $$
\sum_{i=1}^k h(t_{i-1})(t_i-t_{i-1})<{\bar d}_M(p,q)+\ve.
 $$

Consider arbitrary pseudoholomorphic curve
$u_t^\g:D_\d\to M$, which satisfies the conditions
$u_t^\g(0)=\g(t)$ and $(u_t^\g)_*e=\dot\g(t)$.
Define for small $\dd t\in\R_+\subset\C$ the curve
$\hat\g(t;\dd t)= u_t^\g(\dd t)$. Since
$\hat\g(t;\dd t)=\g(t+\dd t)+O(|\dd t|^2)$, propositions~8
and~6 imply that for small $\dd t$ it holds:
  \begin{eqnarray*}
d_M(\g(t),\g(t+\dd t))
&\!\!\!\le\!\!\!&
d_M(\g(t),\hat\g(t;\dd t)) + d_M(\hat\g(t;\dd t), \g(t+\dd t)) \\
&\!\!\!\le\!\!\!&
F_M(\dot\g(t))\dd t + O(|\dd t|^2)\le (1+\ve)h(t)\dd t.
  \end{eqnarray*}

Thus for sufficiently dense partition
  \begin{eqnarray*}
d_M(p,q)\le \sum_{i=1}^k d_M(\g(t_{i-1}),\g(t_i)) <
(1+\ve)({\bar d}_M(p,q)+\ve).
  \end{eqnarray*}
Since $\ve>0$ is arbitrary constant, the theorem is proved. \qed
 \end{proof}

%%%%%%%%%%%%%%%%%%%%%%%%%%%%%%%%%%%%%%%%%%%%%%%%%%%%%%%%%%%%%%%%%%%%%%%%%%%%
%2.3%
\subsection{\hskip-16.5pt . \hskip2pt Hyperbolicity and nonhyperbolicity}

 \begin{dfn}
Almost complex manifold $(M,J)$ is called
{\it hyperbolic} if the pseudodistance $d_M$ is a distance.
 \end{dfn}

Let us consider the unit tangent vectors bundle
$\tau^{(1)}_M:T_1M\to M$ for some norm $|\cdot|$, and let
$F_M^{(1)}:T_1M\to\R$ be the restriction of the Kobayashi-Royden
pseudonorm to it. Proposition~6(i) and theorem~3 imply

 \begin{th}\po
The function $F_M^{(1)}$ is bounded on compact subsets
in $M$. Manifold $M$ is hyperbolic iff $F_M^{(1)}$ is bounded
away from zero on compact subsets.
 \end{th}

Now let us consider the case of nonhyperbolic manifold $M$, for
example let it possess pseudoholomorphic spheres.
In the case of general position for the almost complex structure
$J$, which is tamed by some symplectic form $\om$ on $M$,
the set of all pseudoholomorphic spheres in a fixed homology class
$A\in H_2(M;\Z)$
(completed for compactness by the set of decomposable rational curves)
is a finite-dimensional manifold ${\cal M}(A;J)$ \cite{9, 10}.
We define by the {\it reduction procedure\/} some pseudodistance on
this manifold. Namely for any two pseudoholomorphic spheres
$f_i:S^2\to M$,
defined up to holomorphical reparametrization of $S^2$ let
 $$
d_{{\cal M}}([f_1],[f_2])=d_M(p_1,p_2),
 $$
where $p_i\in\op{Im}(f_i)$ are arbitrary points on the images.
It is easily seen that $d_{{\cal M}}$ is correctly defined
pseudodistance on the manifold $\cal{M}$.

As an example note that the defined pseudodistance $d_{\cal{M}}$
is a distance for almost complex manifold $M^4=\Sigma^2_g\times S^2$
with $g>1$, where the structure $J$ is tamed by the standard
product symplectic form: as in proposition~4 one proves that $M^4$
is fibered by pseudoholomorphic spheres and there is an
isomorphism ${\cal M}\simeq\Sigma^2_g$. However in the case of
four-dimensional manifolds this definition is of importance only
in the case of zero self-intersection. Actually if
$A\cdot A>0$ (for nonexceptional case $A\cdot A\ge0$ \cite{11}),
then two spheres $\op{Im}(f_1)$ and $\op{Im}(f_2)$ of the given
homology class do intersect. Thus $d_{{\cal M}}([f_1],[f_2])=0$.

It was shown in the paper~\cite{16} that for $N$ large enough the
manifold ${\cal M}\times\R^{2N}$ possesses a homotopically canonical
almost complex structure $\td{J}$. Kobayashi pseudodistance
$d_{{\cal M}\times\R^{2N}}$ induces a pseudodistance
$\hat d_{{\cal M}}$ on $\cal{M}$ via reduction over $\R^{2N}$.
In this connection there arises a natural question of existence
of almost complex structures $\td{J}$ such that the
pseudodistances $\hat d_{{\cal M}}$ and $d_{{\cal M}}$ coincide.

%%%%%%%%%%%%%%%%%%%%%%%%%%%%%%%%%%%%%%%%%%%%%%%%%%%%%%%%%%%%%%%%%%%%%%%%%%%%

\ {\hbox to 12.5cm{ \hrulefill }}

{\footnotesize
MSTU n.a. Baumann, Moscow;\quad
kruglikov\verb"@"math.uit.no
}

\end{document}